\documentclass[twoside,10pt]{article}
\usepackage[]{amsmath,amsthm,amssymb}
\usepackage[latin1]{inputenc}

\begin{document}

\begin{center}{{\Large\bf  A class of index transforms generated by} }\end{center}
\begin{center}{{\Large\bf  the Mellin and Laplace operators} }\end{center}
\vspace{0,5cm}

\begin{center}{Semyon  YAKUBOVICH}\end{center}

\markboth{\rm \centerline{ Semyon  YAKUBOVICH}}{}
\markright{\rm \centerline{INDEX TRANSFORMS}}

\begin{abstract} {\noindent Classical integral representation of the
Mellin type kernel
$$ x^{-z}= {1\over \Gamma(z)}\int_0^\infty e^{-xt} t^{z-1}dt,\ x >0,\
{\rm Re\ z} >0,$$ in terms of the Laplace integral gives an idea to
construct a class of non-convolution (index) transforms with the
kernel
$$k^{\pm}_z(x)=\int_0^\infty {e^{-xt^{\pm 1}}\over r(t)} t^{z-1} dt,\ x >0,$$
where $r(t)\neq 0,\ t \in \mathbb{R}_+$ admits a power series
expansion, which has an infinite radius of convergence and the
integral converges absolutely in a half-plane of the complex plane
$z$. Particular examples give the Kontorovich-Lebedev-like
transformation and new transformations with hypergeometric functions
as kernels. Mapping properties and inversion formulas are obtained.
Finally we prove a new inversion theorem for the modified
Kontorovich-Lebedev transform.}

\end{abstract}
\vspace{4mm}

\vspace{2mm}


\vspace{4mm}

\section {Introduction and preliminary results}

In this paper we construct a class of integral transformations of
the non-convolution type, which involves an integration with respect
to parameters of hypergeometric functions. We will base on mapping
and inversion properties of the Mellin and Laplace transforms [5]
given, respectively,  by formulas
$$({\cal M} f)(z)= \int_0^\infty f(t) t^{z-1}dt, \ z \in
\mathbb{C},\eqno(1.1)$$
$$(Lf)(x)= \int_0^\infty f(t)e^{-xt} dt, \ x \in
\mathbb{R}_+,\eqno(1.2)$$ where the integrals converge in an
appropriate sense, which will be clarified below.

The idea to obtain such a new class of index transformations comes
from classical representation of the Mellin kernel
$$x^{-z}= {1\over \Gamma(z)}\int_0^\infty e^{-xt} t^{z-1}dt,\ x >0,\
{\rm Re\ z} >0,\eqno(1.3)$$ where $\Gamma(z)$ is Euler's
gamma-function. Hence doing our steps formally first, we will show
how to invert using (1.3) the modified Mellin transform
$$(Ff)(z)= \int_0^\infty t^{-z} f(t) e^{-at} dt,\ a >0\eqno(1.4)$$
and then will motivate it rigorously in a special class of
functions. Indeed, substituting (1.3) into (1.4) and changing the
order of integration, we find
$$(Ff)(z)= {1\over \Gamma(z)}\int_0^\infty  x^{z-1} \int_0^\infty e^{-(x+a)t}f(t) dtdx.\eqno(1.5)$$
Hence appealing to the inversion formula of the Mellin transform
[4], [5], we have the equality
$${1\over 2\pi i} \int_{\gamma-i\infty}^{\gamma+i\infty} \Gamma(z)
(Ff)(z) x^{-z} dz = \int_0^\infty e^{-(x+a)t}f(t) dt.$$ Substituting
again (1.3) in the left-hand side of the latter equality, we change
the order of integration to get
$$\int_0^\infty e^{-xt} {1\over 2\pi i} \int_{\gamma-i\infty}^{\gamma+i\infty}
(Ff)(z) t^{z-1} dz dt = \int_0^\infty e^{-(x+a)t}f(t)
dt.\eqno(1.6)$$ Finally canceling the Laplace transform (1.2), we
come out with the inversion formula of the transformation (1.4)

$$f(x)= {e^{ax}\over 2\pi i} \int_{\gamma-i\infty}^{\gamma+i\infty} (Ff)(z)
x^{z-1} dz, \ x >0,\eqno(1.7)$$ which coincides with the classical
inversion formula for the Mellin transform up to a simple change of
variables and functions.

The rigorous proof of these reciprocal formulas can be done in a
special class of functions related to the Mellin transform and its
inversion, which was introduced in [7] (see also in [8]). Indeed, we
have

{\bf Definition 1}. Denote by ${\cal M}^{-1}(L_c)$ the space of
functions $f(x), x \in \mathbb{R}_+$,  representable by inverse
Mellin transform of integrable functions $f^{*}(s) \in L_{1}(c)$ on
the vertical line $c =\{s \in \mathbb{C}: {\rm Re} s= c_0\}$:
$$ f(x) = {1\over 2\pi i} \int_c f^{*}(s)x^{-s}ds.\eqno(1.8)$$

The  space ${\cal M}^{-1}(L_c)$  with  the  usual operations  of
addition   and multiplication by scalar is a linear vector space. If
the norm in ${\cal M}^{-1}(L_c)$ is introduced by the formula
$$ \big\vert\big\vert f \big\vert\big\vert_{{\cal
M}^{-1}(L_c)}= {1\over 2\pi }\int^{+\infty}_{-\infty} |
f^{*}\left(c_0 +it\right)| dt,\eqno(1.9)$$
then it becomes  a Banach space.

 {\bf Definition 2 ([7], [8])}.  Let $c_1, c_2 \in \mathbb{R}$
be such that $2 \hbox{sign}\ c_1 + \hbox{sign}\  c_2 \ge 0$. By
${\cal M}_{c_1,c_2}^{-1}(L_c)$ we denote the space of functions
$f(x), x \in \mathbb{R}_+$, representable in the form (1.8), where
$s^{c_2}e^{\pi c_1|s|} f^*(s) \in L_1(c)$.

It is a Banach space with the norm
$$ \big\vert\big\vert f \big\vert\big\vert_{{\cal
M}_{c_1,c_2}^{-1}(L_c)}= {1\over 2\pi }\int_{c} e^{\pi c_1|s|}
|s^{c_2} f^{*}(s) ds|.$$
In particular, letting $c_1=c_2=0$ we get the space ${\cal
M}^{-1}(L_c)$. Moreover, it is easily seen the inclusion $(c_0\neq 0)$

$${\cal M}_{d_1,d_2}^{-1}(L_c) \subseteq {\cal
M}_{c_1,c_2}^{-1}(L_c)$$ when $2 \hbox{sign}(d_1- c_1) + \hbox{sign}
(d_2-c_2) \ge 0$.

We have

{\bf Theorem 1}. {\it Let $f \in {\cal M}^{-1}(L_c),\ a >0, c_0 <
1$. Then transformation $(1.4)$ is well-defined and $(Ff)(z)$ is
analytic in the half-plane ${\rm Re} z <1- c_0$. Moreover,
$$(Ff)(z)={1\over 2\pi i} \int_{c_0-i\infty}^{c_0 +i\infty}
\Gamma(1-s-z)f^*(s)a^{z+s-1} ds,\eqno(1.10)$$
and the operator  $F: {\cal M}^{-1}(L_c) \to L_1({\rm Re} z-i\infty,
{\rm Re} z +i\infty),\ {\rm Re z} < 1-c_0$ is bounded with the norm
satisfying the estimate
$$||F|| \le  a^{{\rm Re} z +c_0-1}\int_{-\infty}^\infty |\Gamma(1-c_0- {\rm Re z}-i\tau)|d\tau.$$
Finally, for all $x >0$ inversion formula $(1.7)$ holds}.
\begin{proof}
Indeed, substituting (1.3) into (1.4) we change the order of
integration by Fubini's theorem via the estimate (see (1.8))

$$|(Ff)(z)|\le \int_0^\infty |t^{-z} f(t)| e^{-at} dt \le {a^{{\rm Re} z +c_0-1}\over 2\pi}\Gamma(1-c_0-{\rm Re} z)
 \int_{c_0-i\infty}^{c_0 +i\infty}|f^*(s)ds| < \infty, $$
which also guarantees the analyticity of $(Ff)(z)$ in the strip
${\rm Re} z <1- c_0$. Thus calculating the inner integral with
respect to $t$ we arrive at the representation (1.10). Finally, the
norm estimation is given by the inequality

$$||Ff||_1= \int_{-\infty}^{\infty} |(Ff)({\rm Re} z + i\tau)|d\tau $$
$$\le {a^{{\rm Re} z +c_0-1}\over 2\pi} \int_{-\infty}^{\infty}
\int_{-\infty}^\infty |\Gamma(1- c_0- {\rm Re z}-i(\tau
+t))f^*(c_0+it)| dtd\tau$$
$$\le a^{{\rm Re} z +c_0-1}\ ||f||_{{\cal M}^{-1}(L_c)}
\int_{-\infty}^\infty |\Gamma(1- c_0- {\rm Re z}- i\tau)|d\tau.$$
In order to prove formula (1.7), we   multiply both sides of (1.10)
by $x^{z-1},\ x>0$ and integrate with respect to $z$ over the line
$(\gamma-i\infty, \gamma +i\infty), \ \gamma < 1-c_0$. Hence
changing the order of integration via the absolute convergence and
calculating the inner integral as an inverse Mellin transform of the
gamma-function, we derive
$$\int_{\gamma-i\infty}^{\gamma+i\infty} (Ff)(z) x^{z-1} dz
={1\over 2\pi i} \int_{c_0-i\infty}^{c_0 +i\infty} f^*(s)a^{s}
\int_{\gamma-i\infty}^{\gamma +i\infty} \Gamma(1-s-z) (ax)^{z-1} dz
ds$$
$$= e^{-ax} \int_{c_0-i\infty}^{c_0 +i\infty} f^*(s) x^{-s}ds= 2\pi i \  e^{-ax}f(x),$$
which proves (1.7).
\end{proof}

\section{General non-convolution transforms}

Let us consider a general non-convolution transformation
$$(Ff)(z)=  \int_0^\infty  k^-_z(x)f(x)dx,\quad  z \in \mathbb{C},\eqno(2.1)$$
where
$$k^{-}_z(x)=\int_0^\infty {e^{-x/t}\over r(t)} t^{z-1} dt,\ x >0\eqno(2.2)$$
and $r(t)\neq 0, t >0$ admits the series representation $r(t)=
\sum_{k=0}^\infty a_kt^k$ with an infinite radius of convergence.

{\bf Theorem 2}. {\it Let $f \in {\cal M}^{-1}(L_c),\ c_0 < 1$. Let
$r^{-1}(t) \in L_1(\mathbb{R}_+; t^{\gamma- c_0}dt), \gamma \in
\mathbb{R}, \ \rho(s) \in L_1(1+\gamma-c_0 -i\infty, 1+\gamma-c_0
+i\infty),$ where $\rho(s)$ is the Mellin transform $(1.1)$ of the
function $r^{-1}(t)$
$$\rho(s)= \int_0^\infty {t^{s-1}\over r(t)} dt.\eqno(2.3)$$
Then transformation $(2.1)$ is well-defined and $(Ff)(z), {\rm Re}
z=\gamma$ can be represented in the form
$$(Ff)(z)={1\over 2\pi i} \int_{c_0-i\infty}^{c_0 +i\infty}
\Gamma(1-s)\rho(1+z-s) f^*(s) ds.\eqno(2.4)$$ The operator $(2.1)$
is bounded from ${\cal M}^{-1}(L_c)$ into $L_1(1-c_0+\gamma-i\infty,
1-c_0+\gamma +i\infty)$ and
$$||F|| \le  \Gamma(1-c_0)\int_{-\infty}^\infty |\rho(1-c_0 +\gamma -i\tau)|d\tau.\eqno(2.5)$$
Moreover, when $\gamma < 0$ and $(Ff)(z)/\Gamma(-z) \in
L_1(\gamma-i\infty, \gamma+ i\infty)$, for all $x
>0$ the inversion formula holds
$$f(x)= {1\over 2\pi i} \int_{\gamma-i\infty}^{\gamma +i\infty}
\hat{k}_z^-(x)(Ff)(z) dz,\eqno(2.6)$$
where $$\hat{k}_z^-(x)= \sum_{k=0}^\infty  {a_k x^{k-z-1}\over
\Gamma(k-z)}\eqno(2.7)$$ and integral $(2.6)$ converges absolutely}.
\begin{proof}
Since using (1.8), (2.2) and conditions of the theorem
$$|(Ff)(z)| \le \int_0^\infty  |k^-_z(x)f(x)|dx \le {1\over
2\pi}\int_0^\infty |k^-_z(x)| x^{-c_0}dx
 \int_{c_0-i\infty}^{c_0 +i\infty}|f^*(s)ds|$$
$$\le {1\over
2\pi}\int_0^\infty \int_0^\infty {e^{-x/t} x^{-c_0}\over |r(t)|}
t^{{\rm Re} z-1} dtdx
 \int_{c_0-i\infty}^{c_0 +i\infty}|f^*(s)ds|$$
$$= {\Gamma(1-c_0)\over 2\pi}\int_0^\infty {t^{{\rm Re} z-c_0}\over |r(t)|}
 dt  \int_{c_0-i\infty}^{c_0 +i\infty}|f^*(s)ds| < \infty,$$
one can substitute (1.8) into (2.1) and change the order of
integration via the absolute convergence. After calculation of the
inner integral employing the convolution property of the Mellin
transform [5] and minding (2.3), we come out with  representation
(2.4). Hence
$$||Ff||_1 \le {\Gamma(1-c_0)\over 2\pi} \int_{c_0-i\infty}^{c_0 +i\infty}
\int_{-\infty}^\infty \left|\rho(1+\gamma-c_0-i\tau) f^*(s)
ds\right|d\tau$$
$$\le \Gamma(1-c_0)||f||_{{\cal M}^{-1}(L_c)} \int_{-\infty}^\infty
\left|\rho(1+\gamma-c_0-i\tau)\right|d\tau,$$
which yields (2.5). Returning to (2.4) we multiply both sides of
this equality by $x^z,\ x>0$ and integrate with respect to $z$ over
the line $(\gamma-i\infty, \gamma +i\infty)$. Changing the order of
integration by Fubini's theorem, which is applicable by virtue of
the absolute convergence of the corresponding integral, we calculate
the inner integral via the inversion theorem for the Mellin
transform [5], since the original function and its image are
integrable. Therefore reciprocally from (2.3) for all $x>0$ we have
$${1\over 2\pi i} \int_{\gamma-i\infty}^{\gamma +i\infty}
\rho(1+z-s) x^z dz= x^{s-1} [r(1/x)]^{-1}$$
and
$$\int_{\gamma-i\infty}^{\gamma +i\infty} (Ff)(z) x^z dz = [r(1/x)]^{-1} \int_{c_0-i\infty}^{c_0 +i\infty}
\Gamma(1-s)f^*(s)x^{s-1} ds.$$
Hence appealing in the right-hand side of the latter equality to the
Parseval identity for the Mellin transform, we derive
$${r(1/x)\over 2\pi i} \int_{\gamma-i\infty}^{\gamma +i\infty} (Ff)(z) x^z
dz = \int_0^\infty e^{-xt} f(t)dt.\eqno(2.8)$$
Meanwhile, bearing in mind the series representation of the function
$r(1/x)$ and its infinite radius of convergence, we observe that the
series of coefficients  $\sum_{k=0}^\infty a_k$ converges
absolutely. So, the left-hand side of (2.8) can be treated as
follows
$${1\over 2\pi i} \int_{\gamma-i\infty}^{\gamma +i\infty} (Ff)(z) \sum_{k=0}^\infty  a_k
x^{z-k} dz = {1\over 2\pi i} \int_{\gamma-i\infty}^{\gamma +i\infty}
(Ff)(z) \sum_{k=0}^\infty  {a_k \over \Gamma(k-z)} \int_0^\infty
e^{-xt}t^{k-z-1}dt dz$$
$$= \int_0^\infty e^{-xt} {1\over 2\pi i} \int_{\gamma-i\infty}^{\gamma +i\infty}
(Ff)(z)\  \sum_{k=0}^\infty  {a_k t^{k-z-1}\over \Gamma(k-z)} \
dzdt,$$
where the change of the order of integration is possible via the
integrability of the function $(Ff)(z)/\Gamma(-z)$ and the estimate
$(\gamma < 0)$
$$\int_{\gamma-i\infty}^{\gamma +i\infty}\left|(Ff)(z) \sum_{k=0}^\infty  {a_k \over \Gamma(k-z)} \int_0^\infty
e^{-xt}t^{k-z-1}dt dz\right| \le \int_{\gamma-i\infty}^{\gamma
+i\infty} \left|{(Ff)(z) \over \Gamma(-z)}dz\right|\sum_{k=0}^\infty
{|a_k| \Gamma(k-\gamma) \over (-\gamma)_k} $$
$$= \Gamma(-\gamma) \int_{\gamma-i\infty}^{\gamma
+i\infty} \left|{(Ff)(z) \over \Gamma(-z)}dz\right|
\sum_{k=0}^\infty |a_k| < \infty, $$ where $(-\gamma)_k$ is the
Pochhammer symbol [1], Vol. I.

Thus returning to (2.8), we get the equality
$$\int_0^\infty e^{-xt} {1\over 2\pi i} \int_{\gamma-i\infty}^{\gamma
+i\infty} (Ff)(z)\  \sum_{k=0}^\infty  {a_k t^{k-z-1}\over
\Gamma(k-z)} \ dzdt =\int_0^\infty e^{-xt} f(t)dt, \ x
>0.\eqno(2.9)$$
As we see, equality (2.9) is true for all $x >0$, where functions
under the convergent Laplace integrals in its both sides are
continuous on $\mathbb{R}_+$ owing to  the condition $f \in {\cal
M}^{-1}(L_c)$ and assumptions of the theorem. Therefore one can
cancel the Laplace transform (1.2) in (2.9) by virtue of the
uniqueness theorem (see in [2]). Consequently, we established  the
inversion formula (2.6) and completed the proof of Theorem 2.
\end{proof}

We note, that if $r(t)\equiv P_n(t)$ is a polynomial with no zeros
on the line $\mathbb{R}_+$, then Theorem 2 is true, where instead of
the series final sums are involved. But in this case we can also
prove similarly the theorem about mapping properties and inversion
formula of the following transformation
$$(Ff)_n(z)=  \int_0^\infty  k_n(z,x)f(x)dx,\quad  z \in \mathbb{C},\ n \in \mathbb{N},\eqno(2.10)$$
where
$$k_n(z,x)=\int_0^\infty {e^{-xt}\over P_n(t)} t^{z-1} dt,\ x >0\eqno(2.11)$$
and $P_n(t)= \sum_{k=0}^n a_k t^k$.  Precisely, we have the
following result.

{\bf Theorem 3.} {\it Let $n \in \mathbb{N}, \ f \in {\cal
M}^{-1}(L_c),\ c_0 < 1$. Let $\gamma \in \mathbb{R}, -c_0< \gamma <
n-c_0, \  \rho_n(s) \in L_1(\gamma+c_0 -i\infty, \gamma+c_0
+i\infty),$ where $\rho_n(s)$ is the Mellin transform $(1.1)$ of
$P_n^{-1}(t)$
$$\rho_n(s)= \int_0^\infty {t^{s-1}\over P_n(t)} dt.\eqno(2.12)$$
Then transformation $(2.10)$ is well-defined and $(Ff)_n(z)$ has the
representation
$$(Ff)_n(z)={1\over 2\pi i} \int_{c_0-i\infty}^{c_0 +i\infty}
\Gamma(1-s)\rho_n(s+z) f^*(s) ds.\eqno(2.13)$$ The operator $(2.10)$
is bounded from ${\cal M}^{-1}(L_c)$ into $L_1(c_0+\gamma-i\infty,
c_0+\gamma +i\infty)$ and
$$||F|| \le  \Gamma(1-c_0)\int_{-\infty}^\infty |\rho_n(c_0 +\gamma -i\tau)|d\tau.$$
Moreover, when $\hbox{max}\  (-c_0, \ n-1) < \gamma < n-c_0$ and
$(Ff)_n(z)/\Gamma(1+z) \in L_1(\gamma-i\infty, \gamma+ i\infty)$,
for all $x
>0$ the inversion formula holds
$$f(x)= {1\over 2\pi i} \int_{\gamma-i\infty}^{\gamma +i\infty}
\hat{k}_n(z,x)(Ff)(z) dz,\eqno(2.14)$$
where $$\hat{k}_n(z,x)= \sum_{k=0}^n  {a_k x^{z-k}\over
\Gamma(1+z-k)}\eqno(2.15)$$ and integral $(2.14)$ converges
absolutely}.

\section{Examples of new index transforms}

We start this section showing an interesting example of the
transform (2.1) recently discovered by the author (see in [10]). In
fact, let $r(t)= e^{-t}$. Then calculating the integral (2.2), we
get $k^{-}_z(x)= 2x^{z/2}K_z(2\sqrt x)$, where $K_z(2\sqrt x)$ is
the modified Bessel function [1], Vol. II.  As it is easily seen,
integral (2.1) converges absolutely for any $x \in \mathbb{R}_+, z
\in \mathbb{C}$ and represents an entire function of $z$. On the
other hand, the kernel $k^{-}_z(x)$ can be written with the use of
the Parseval relation for the Mellin transform, which leads to the
representation
$$ 2 x^{z/2}K_z(2\sqrt x)= {1\over 2\pi i}
\int_{\gamma-i\infty}^{\gamma+i\infty} \Gamma(s+z)\Gamma(s) x^{-s}
ds,\ x > 0.\eqno(3.1)$$
Thus the  transformation (2.1) in this case has the form
$$(Ff)(z)=  2\int_0^\infty  x^{z/2}K_z(2\sqrt x)f(x)dx.\eqno(3.2)$$
This transformation looks like the Kontorovich-Lebedev transform
[4], [8], [9]. However, it is a completely different operator and
cannot be reduced to the Kontorovich-Lebedev integral by any change
of variables and functions. As far as the author is aware, the
transform (3.2) was not studied yet, taking into account his mapping
properties and inversion formula in an appropriate class of
functions. An analog of Theorem 2 for this case is

{\bf Theorem 4}\ [10]. {\it Let $f \in {\cal M}^{-1}(L_c)$ and $c_0
< 1$. Then transformation $(3.2)$ is well-defined and $(Ff)(z)$ is
analytic in the half-plane ${\rm Re} z > c_0 -1$. Further,
$$(Ff)(z)={1\over 2\pi i}
\int_{c_0-i\infty}^{c_0 +i\infty} \Gamma(1-s+z)\Gamma(1-s)f^*(s)
ds,\eqno(3.3)$$
and the operator  $F: {\cal M}^{-1}(L_c) \to L_1(\gamma -i\infty,
\gamma +i\infty),\ \gamma > c_0 - 1$ is bounded with the norm
satisfying the estimate
$$||F|| \le  \Gamma(1-c_0) \int_{-\infty}^\infty |\Gamma(1-c_0+ \gamma + i\tau)|d\tau.$$

Moreover, when $c_0-1 < \gamma < 0$ and $(Ff)(z)/\Gamma(-z) \in
L_1(\gamma-i\infty, \gamma+ i\infty)$, for all $x
>0$ the inversion formula holds
$$f(t)= {1\over 2\pi i} \int_{\gamma-i\infty}^{\gamma +i\infty} I_{-(1+z)}\left(2\sqrt t\right)  t^{-(1+z)/2}\ (Ff)(z)
dz,\eqno(3.4)$$
where $I_\nu(w)$ is the modified Bessel function of the third kind
[2], Vol. II and the integral $(3.4)$ converges absolutely.}

The Kontorovich-Lebedev-like transformation (3.2) can be generalized
considering the following kernel
$$S_m(x,z)= \int_0^\infty e^{-{x\over t}- t^m}
t^{z-1}dt, \ x >0, \ m \in \mathbb{N}.\eqno(3.5)$$
We calculate integral (3.5) in terms of the Meijer $G$-function [1],
Vol. I. Precisely, we derive
$$\int_0^\infty e^{-{x\over t}- t^m} t^{z-1}dt ={1\over 2\pi i m}
\int_{\nu-i\infty}^{\nu +i\infty} \Gamma\left({s+z\over
m}\right)\Gamma(s) x^{-s} ds.$$ Appealing to the Gauss-Legendre
multiplication formula for gamma-function [1], Vol.I
$$\Gamma(ms)= m^{ms-1/2}(2\pi)^{(1-m)/2}\prod_{k=0}^{m-1}
\Gamma\left(s+ {k\over m}\right), \ m \in \mathbb{N},$$ the latter
Mellin-Barnes integral becomes the following Meijer $G$-function
$${1\over 2\pi i m} \int_{\nu-i\infty}^{\nu +i\infty}
\Gamma\left({s+z\over m}\right)\Gamma(s) x^{-s} ds$$ $$ =
{(2\pi)^{(1-m)/2}\over 2\pi i\  m^{1/2}} \int_{{\nu\over
m}-i\infty}^{{\nu\over m} +i\infty} \Gamma\left(s+{z\over
m}\right)\prod_{k=0}^{m-1} \Gamma\left(s+{k\over m}\right)
\left({x\over m}\right)^{-ms} ds$$
$$=  {(2\pi)^{(1-m)/2}\over  m^{1/2}}G{m+1,0 \atop 0,m+1}
 \left( \left({x\over m}\right)^{m} \bigg\vert { - \atop
0, {1\over m},\dots, {m-1\over m}, {z\over m}} \right).$$ Thus we
come out with the more general index transformation, namely
$$(Ff)(z)= {(2\pi)^{(1-m)/2}\over  m^{1/2}}\int_0^\infty
G{m+1,0 \atop 0,m+1} \left( \left({x\over m}\right)^{m} \bigg\vert {
- \atop 0, {1\over m},\dots, {m-1\over m}, {z\over m}}
\right)f(x)dx,\ z \in \mathbb{C}.\eqno(3.6)$$
The inversion formula (2.6) for this case is given accordingly
$$f(x)= {1\over 2\pi i} \int_{\gamma-i\infty}^{\gamma +i\infty}
\hat{S}_m(z,x)(Ff)(z) dz,\quad  x >0,\eqno(3.7)$$
where the kernel $\hat{S}_m(z,x)$ can be expressed in terms of the
hyper-Bessel functions. In fact,  the corresponding kernel (2.7) is
written as the generalized hypergeometric series
$$\hat{S}_m(z,x)= \sum_{n=0}^\infty  {x^{mn-z-1}\over n!\Gamma(mn-z)}
= (2\pi)^{(m-1)/2}m^{1/2}x^{-1}\sum_{n=0}^\infty {(xm)^{mn-z} \over
n!\prod_{k=0}^{m-1}\Gamma\left(n+ {k-z\over m}\right)}$$ $$
=\frac{(2\pi)^{(m-1)/2}m^{1/2}x^{-1}}{(xm)^z\prod_{k=0}^{m-1}\Gamma\left(
{k-z\over m}\right)}{}_0F_m\left({-z\over m}, {1-z\over
m},\dots,{m-1-z\over m}; (xm)^m \right).$$

We have

{\bf Theorem 5}. {\it Let $m \in \mathbb{N}, f \in {\cal
M}^{-1}(L_c)$ and $c_0 < 1$. Then transformation $(3.6)$ is
well-defined and $(Ff)(z)$ is analytic in the half-plane ${\rm Re} z
> c_0 -1$. Further,
$$(Ff)(z)={1\over 2\pi i m}
\int_{c_0-i\infty}^{c_0 +i\infty} \Gamma\left({1-s+z\over
m}\right)\Gamma(1-s)f^*(s) ds,$$
and the operator  $F: {\cal M}^{-1}(L_c) \to L_1(\gamma-i\infty,
\gamma +i\infty),\ \gamma > c_0 - 1$ is bounded with the norm
satisfying the estimate
$$||F|| \le  \Gamma(1-c_0) \int_{-\infty}^\infty \left|\Gamma\left({1\over m}(1-c_0+ \gamma + i\tau)\right)\right|d\tau.$$

Moreover, when $c_0-1 < \gamma < 0$ and
$(Ff)(z)/\prod_{k=0}^{m-1}\Gamma\left( {k-z\over m}\right) \in
L_1(\gamma-i\infty, \gamma+ i\infty)$, for all $x
>0$ the inversion formula $(3.7)$ holds with the absolutely convergent integral}.

Next, calling relation (2.3.6.9) in [3], Vol. 1
$$\int_0^\infty {e^{-xt}\over (t+1)^n} t^{z-1}dt=
\Gamma(z)\Psi(z,z+1-n;x),\ x >0,\ {\rm Re} z >0,\ n \in
\mathbb{N},\eqno(3.8)$$
where $\Psi(a,b;w)$ is Tricomi's function [1], Vol. I, we are ready
to introduce the corresponding index transform (2.10) in the form
$$(Ff)(z)= \Gamma(z) \int_0^\infty  \Psi(z,z+1-n;x)f(x)dx.\eqno(3.9)$$
According to Theorem 3 and equalities (2.14), (2.15)  its inversion
formula is given by the following integral
$$f(x)= {1\over 2\pi i} \int_{\gamma-i\infty}^{\gamma +i\infty}\left(\sum_{k=0}^n  {n\choose k}{x^{z-k}\over
\Gamma(1+z-k)}\right) (Ff)(z) dz,\ x >0.\eqno(3.10)$$
But the finite sum inside (3.10) can be expressed in terms of the
generalized Laguerre polynomials. In fact, appealing to relation
(7.17.1.1) in [3], Vol. 3 we obtain
$$\sum_{k=0}^n  {n\choose k}{x^{z-k}\over \Gamma(1+z-k)}= {x^z\over
\Gamma(1+z)}\ {}_2F_0\left(-n,-z; {1\over x}\right)= {x^{z-n}
n!\over \Gamma(1+z)}\ L_n^{z-n}(-x).$$

Hence an analog of Theorem 3 is

{\bf Theorem 6.} {\it Let $n \in \mathbb{N}, \ f \in {\cal
M}^{-1}(L_c),\ c_0 < 1$. Let $\gamma \in \mathbb{R}, 0< \gamma <
n-c_0.$ Then transformation $(3.9)$ is well-defined and $(Ff)(z)$
has the representation
$$(Ff)(z)={1\over 2\pi i(n-1)!} \int_{c_0-i\infty}^{c_0 +i\infty}
\Gamma(s+z)\Gamma(n-s-z)\Gamma(1-s) f^*(s) ds.$$ The operator
$(3.9)$ is bounded from ${\cal M}^{-1}(L_c)$ into
$L_1(c_0+\gamma-i\infty, c_0+\gamma +i\infty)$ and
$$||F|| \le  {\Gamma(1-c_0)\over (n-1)!}\int_{-\infty}^\infty |\Gamma(c_0+\gamma+i\tau)\Gamma(n-c_0-\gamma-i\tau)|d\tau.$$
Moreover, if $ n-1 < \gamma < n-c_0$ and $(Ff)(z)/\Gamma(1+z) \in
L_1(\gamma-i\infty, \gamma+ i\infty)$, for all $x >0$ the inversion
formula $(3.10)$ holds
$$f(x)= {n!\over 2\pi i} \int_{\gamma-i\infty}^{\gamma +i\infty} {x^{z-n}
\over \Gamma(1+z)}\ L_n^{z-n}(-x)(Ff)(z) dz,\ x >0,$$
where the latter integral is absolutely convergent}.

In particular, letting $n=1$ in (3.8), we use relation (2.3.6.13) in
[3], Vol. 1 to get
$$\int_0^\infty {e^{-xt}\over t+1} t^{z-1}dt=
\Gamma(z)e^x\ \Gamma(1-z, x),\ x >0,\ {\rm Re} z
>0,\eqno(3.11)$$
where $\Gamma(w, a)$ is incomplete gamma-function [1], Vol. I.
Consequently, we have the reciprocal pair of index transforms
$$(Ff)(z)= \Gamma(z)\int_0^\infty \Gamma(1-z, x)e^x
f(x)dx,\eqno(3.12)$$
$$f(x)= {1\over 2\pi i} \int_{\gamma-i\infty}^{\gamma +i\infty} {x^{z}\over \Gamma(z)}\left[{1\over z}+ {1\over x}\right]
 (Ff)(z) dz,\ x >0.\eqno(3.13)$$
{\bf Theorem 7.} {\it Let $f \in {\cal M}^{-1}(L_c),\ c_0 < 1$. Let
$\gamma \in \mathbb{R}, 0< \gamma < 1-c_0.$ Then transformation
$(3.12)$ is well-defined and $(Ff)(z)$ has the representation
$$(Ff)(z)={1\over 2i} \int_{c_0-i\infty}^{c_0 +i\infty}
{\Gamma(1-s)\over \sin(\pi(s+z))} f^*(s) ds.$$ The operator $(3.12)$
is bounded from ${\cal M}^{-1}(L_c)$ into $L_1(c_0+\gamma-i\infty,
c_0+\gamma +i\infty)$ and
$$||F|| \le  \pi\Gamma(1-c_0)\int_{-\infty}^\infty {d\tau\over |\sin(\pi(c_0+\gamma+i\tau))|}.$$
Moreover, if $0 < \gamma < 1-c_0$ and $(Ff)(z)/\Gamma(1+z) \in
L_1(\gamma-i\infty, \gamma+ i\infty)$, for all $x >0$ the inversion
formula $(3.13)$ holds, where the integral is absolutely
convergent}.

A final example, which we are going to consider is generated by
relation (2.3.7.8) in [3], Vol. 1. It involves a combination of
hypergeometric functions ${}_1F_2$, namely
$$H(z,x)= \int_0^\infty {e^{-xt}\over (t^2+1)^n}t^{z-1}dt
= x^{2n-z}\Gamma(z- 2n)\ {}_1F_2\left(n; 1+n-{z\over 2}, n+
{1-z\over 2}; - {x^2\over 4}\right)$$ $$ + \frac{\Gamma\left(n-
{z\over 2}\right)\Gamma\left({z\over 2}\right)}{2(n-1)!}\
{}_1F_2\left({z\over 2}; {1\over 2},  1- n+{z\over 2}; - {x^2\over
4}\right)$$$$- \frac{x\Gamma\left(n- {z+1\over
2}\right)\Gamma\left({z+1\over 2}\right)}{2(n-1)!}\
{}_1F_2\left({z+1\over 2}; {3\over 2},  {3+z\over 2}-n; - {x^2\over
4}\right),\ n \in \mathbb{N}, \ x >0, \ {\rm Re} z > 0.\eqno(3.14)$$
So, the index transformation with the kernel (3.14)
$$(Ff)(z)=  \int_0^\infty  H(z,x)f(x)dx\eqno(3.15)$$
admits the following inversion formula
$$f(x)= {1\over 2\pi i} \int_{\gamma-i\infty}^{\gamma +i\infty}\left(\sum_{k=0}^n  {n\choose k}{x^{z-2k}\over
\Gamma(1+z-2k)}\right) (Ff)(z) dz,\ x >0,\eqno(3.16)$$
In the meantime,
$$\sum_{k=0}^n  {n\choose k}{x^{z-2k}\over \Gamma(1+z-2k)}= {x^z\over
\Gamma(1+z)}\ {}_3F_0\left(-n,-{z\over 2}, {1-z\over 2}; - {4\over
x^2}\right).$$

Therefore we arrive at the following result.

{\bf Theorem 8.} {\it Let $n \in \mathbb{N}, \ f \in {\cal
M}^{-1}(L_c),\ c_0 < 1$. Let $\gamma \in \mathbb{R}, 0< \gamma <
n-c_0.$ Then transformation $(3.15)$ is well-defined and $(Ff)(z)$
has the representation
$$(Ff)(z)={1\over 4\pi i(n-1)!} \int_{c_0-i\infty}^{c_0 +i\infty}
\Gamma\left({s+z\over 2}\right)\Gamma\left(n- {s+z\over
2}\right)\Gamma(1-s) f^*(s) ds.$$ The operator $(3.15)$ is bounded
from ${\cal M}^{-1}(L_c)$ into $L_1(c_0+\gamma-i\infty, c_0+\gamma
+i\infty)$ and
$$||F|| \le  {\Gamma(1-c_0)\over 2(n-1)!}\int_{-\infty}^\infty
\left|\Gamma\left({c_0+\gamma+i\tau\over 2}\right)\Gamma\left(n-
{c_0+\gamma+i\tau\over 2}\right)\right|d\tau.$$
Moreover, if $n-1 < \gamma < n-c_0$ and $(Ff)(z)/\Gamma(1+z) \in
L_1(\gamma-i\infty, \gamma+ i\infty)$, for all $x >0$ the inversion
formula $(3.16)$ holds
$$f(x)= {1\over 2\pi i} \int_{\gamma-i\infty}^{\gamma +i\infty} {x^z\over
\Gamma(1+z)}\ {}_3F_0\left(-n,-{z\over 2}, {1-z\over 2}; - {4\over
x^2}\right)(Ff)(z) dz,\ x >0,$$
where the latter integral is absolutely convergent}.

\section{A new inversion theorem for the modified
Kontorovich-Lebedev transform}

In this final section we will prove a new inversion theorem for the
modified Kontorovich-Lebedev transform
$$g(x)= \int_0^\infty e^{-x/2}K_{i\tau}\left({x\over
2}\right)f(x)dx,\ \tau \in \mathbb{R}.\eqno(4.1)$$
It is easily seen that making elementary changes of variables and
functions we come out with the classical Kontorovich-Lebedev
transform (cf. [4], [8], [9]). We will prove that an arbitrary
function $f \in {\cal M}_{1/2,\nu}^{-1}(L_c),\ \nu \in \mathbb{R}$
(see Definition 2)  can be expanded for all $x>0$ in terms of the
following iterated Kontorovich-Lebedev integral
$$f(x)= {e^{x/2}\over \pi^2 x}\int_{-\infty}^\infty \tau\sinh\pi\tau K_{i\tau}\left({x\over 2}\right)\int_0^\infty e^{-y/2}K_{i\tau}\left({y\over
2}\right)f(y)dy.\eqno(4.2)$$
We note,  that similar expansion in the space ${\cal
M}_{0,1/4}^{-1}(L)$ was studied in [6].

{\bf Theorem 9}. {\it Let $f \in {\cal M}_{1/2,\nu}^{-1}(L_c),\ \nu
> c_0-1,\ 0< c_0 <1$. The expansion $(4.2)$ holds for all $x>0$, where
the integral with respect to $y$ is absolutely convergent and the
integral with respect to $\tau$ exists in the Riemann improper
sense.}

\begin{proof} Calling relation (8.4.23.3) in [3], Vol. 3, we use the
Parseval equality for the Mellin transform [5] and Definition 2 to
write  integral (4.1) in the form
$$\int_0^\infty e^{-y/2}K_{i\tau}\left({y\over
2}\right)f(y)dy = {\sqrt \pi\over 2\pi
i}\int_{c_0-i\infty}^{c_0+i\infty}
\frac{\Gamma(1-s-i\tau)\Gamma(1-s+i\tau)}{\Gamma(3/2-s)}f^*(s)ds.\eqno(4.3)$$
But due to Stirling's formula for gamma-functions [1], Vol. I, we
have $\Gamma(3/2-s)= O(|s|^{1-c_0}e^{-\pi|s|/2}), \ |s| \to \infty$.
Therefore via conditions of the theorem the function
$$\frac{f^*(s)}{\Gamma(3/2-s)}$$
is Lebesgue integrable over the line $(c_0- i\infty, c_0+i\infty)$.
Hence denoting by
$$h(x)= {1\over 2\pi i}\int_{c_0-i\infty}^{c_0+i\infty}
\frac{f^*(s)}{\Gamma(3/2-s)} x^{-s}ds\eqno(4.4)$$
and applying again the Parseval equality for the Mellin transform
together with representation (3.1) for the modified Bessel function,
the right-hand side of (4.3) becomes
$${\sqrt \pi\over 2\pi i}\int_{c_0-i\infty}^{c_0+i\infty}
\frac{\Gamma(1-s-i\tau)\Gamma(1-s+i\tau)}{\Gamma(3/2-s)}f^*(s)ds=
2\sqrt\pi \int_0^\infty K_{2i\tau}(2\sqrt y)h(y)dy.\eqno(4.5)$$
Substituting the right-hand side of (4.5) into (4.2), we change the
order of integration by Fubini's theorem, which is applicable owing
to the absolute convergence of the iterated integral. Indeed, fixing
a positive $x$, we appeal to the uniform inequality for the modified
Bessel function [9]
$$|K_{i\tau}(x)| \le e^{-\delta|\tau|}K_0(x\cos\delta), \ \delta \in
\left[0, {\pi\over 2}\right[,$$
definition (4.4) of $h(x)$, and asymptotic behavior of the modified
Bessel function [1], Vol. II   to have the estimate
$$\int_{-\infty}^\infty \left|\tau\sinh\pi\tau K_{i\tau}\left({x\over 2}\right)\right|
\int_0^\infty \left|K_{2i\tau}(2\sqrt y)h(y)\right|dy d\tau \le
{1\over 2\pi} K_{0}\left({x\cos\delta\over
2}\right)\int_{-\infty}^\infty |\tau\sinh\pi\tau |
e^{-3\delta|\tau|} d\tau $$ $$\times \int_0^\infty K_{0}(2\cos\delta
\sqrt y)y^{-c_0}dy \int_{c_0-i\infty}^{c_0+i\infty}
\left|\frac{f^*(s)}{\Gamma(3/2-s)} ds\right| < \infty,$$
since one can choose $\delta \in ]\pi/3, \pi/2[.$ Calculating the
inner index integral via relation (2.16.52.9) in [3], Vol. 2, which
is slightly corrected by the author
$$\int_{-\infty}^\infty \tau\sinh\pi\tau K_{i\tau}\left({x\over 2}\right)K_{2i\tau}(2\sqrt y) d\tau
= {1\over 2} \sqrt{{\pi^3 y\over x}}\ e^{-x/2- y/x},$$
we take the result, writing the right-hand side of (4.2) in the form
$${1\over x^{3/2}} \int_0^\infty e^{-y/x} h(y)\sqrt y \ dy.$$
Meanwhile, using expression (4.4) for $h(x)$ and changing the order
of integration after its substitution in the latter integral, we
easily deduce the equalities
$${1\over x^{3/2}} \int_0^\infty e^{-y/x} h(y)\sqrt y \ dy=
{1\over 2\pi i \ x^{3/2}}\int_{c_0-i\infty}^{c_0+i\infty}
\frac{f^*(s)}{\Gamma(3/2-s)} \int_0^\infty e^{-y/x} y^{1/2-s}dy ds$$
$$={1\over 2\pi i}\int_{c_0-i\infty}^{c_0+i\infty}
f^*(s) x^{-s}ds = f(x)$$
via Definition 2. Thus we proved (4.2) and completed the proof of
the theorem.
\end{proof}

\bigskip
\centerline{{\bf Acknowledgments}}
\bigskip
The present investigation was supported, in part,  by the "Centro de
Matem{\'a}tica" of the University of Porto.

\bigskip
\centerline{{\bf References}}
\bigskip
\baselineskip=12pt
\medskip
\begin{enumerate}

\item[{\bf 1.}\ ]
 A. Erd\'elyi, W. Magnus, F. Oberhettinger and F.G. Tricomi,
{\it Higher Transcendental Functions}, Vols. I and  II, McGraw-Hill,
New York, London and Toronto (1953).

\item[{\bf 2.}\ ] V.A. Ditkin and A.P. Prudnikov, {\it Operational
Calculus}. Nauka, Moscow, 1975 (in Russian).

\item[{\bf 3.}\ ]  A.P. Prudnikov, Yu. A. Brychkov and O. I.
Marichev, {\it Integrals and Series: Vol. 1: Elementary Functions},
Gordon and Breach, New York (1986); {\it Vol. 2: Special Functions},
Gordon and Breach, New York (1986); {\it Vol. 3: More Special
Functions}, Gordon and Breach, New York (1990).

\item[{\bf 4.}\ ]I.N. Sneddon, {\it The Use of Integral
Transforms}, McGray Hill, New York (1972).

\item[{\bf 5.}\ ]  E.C. Titchmarsh, {\it  An Introduction to the
Theory of Fourier Integrals}, Clarendon Press, Oxford ( 1937).

\item[{\bf 6.}\ ] Vu Kim Tuan and S.B. Yakubovich,  The Kontorovich-Lebedev integral
transform in a new class of functions, {\it Dokl. Akad. Nauk BSSR},
{\bf 29} (1985), 11-14.   (in Russian);  translation in {\it Amer.
Math. Soc. Transl.}  {\bf 137} (1987), 61-65.

\item[{\bf 7.}\ ] Vu Kim Tuan, O.I. Marichev and S.B. Yakubovich, Composition structure of integral
transformations,  {\it J. Soviet Math.}, {\bf 33} (1986), 166-169.

\item[{\bf 8.}\ ] S. B. Yakubovich and  Yu. F. Luchko, {\it The
Hypergeometric Approach to Integral Transforms and Convolutions}.
Mathematics and its Applications, 287. Kluwer Academic Publishers
Group, Dordrecht (1994).

\item[{\bf 9.}\ ] S.B. Yakubovich, {\it Index Transforms},  World
Scientific Publishing Company, Singapore, New Jersey, London and
Hong Kong (1996).

\item[{\bf 10.}\ ] S.Yakubovich, A new Kontorovich-Lebedev-like transformation,
{\it Commun. Math. Anal.} {\bf 13} (2012), N 1, 86-99.

\end{enumerate}

\vspace{5mm}

\noindent S.Yakubovich\\
Department of Mathematics,\\
Faculty of Sciences,\\
University of Porto,\\
Campo Alegre st., 687\\
4169-007 Porto\\
Portugal\\
E-Mail: syakubov@fc.up.pt\\

\end{document}